\documentclass[12pt]{article}
\usepackage{amsfonts,amssymb,epsfig}
\usepackage{amsmath}

\date{}
\begin{document}
\newtheorem{df}{Definition}
\newtheorem{thm}{Theorem}

\newtheorem{lm}{Lemma}
\newtheorem{pr}{Proposition}
\newtheorem{co}{Corollary}
\newtheorem{re}{Remark}
\newtheorem{note}{Note}
\newtheorem{claim}{Claim}
\newtheorem{problem}{Problem}

\def\R{{\mathbb R}}

\def\E{\mathbb{E}}
\def\calF{{\cal F}}
\def\N{\mathbb{N}}
\def\calN{{\cal N}}
\def\calH{{\cal H}}
\def\n{\nu}
\def\a{\alpha}
\def\d{\delta}
\def\t{\theta}
\def\e{\varepsilon}
\def\t{\theta}
\def\g{\gamma}
\def\G{\Gamma}
\def\b{\beta}
\def\pf{ \noindent {\bf Proof: \  }}

\newcommand{\qed}{\hfill\vrule height6pt
width6pt depth0pt}
\def\endpf{\qed \medskip} \def\colon{{:}\;}
\setcounter{footnote}{0}

\renewcommand{\qed}{\hfill\vrule height6pt  width6pt depth0pt}

\title{A Schauder basis for $L_1(0,\infty)$ consisting of non-negative functions}

\author{William B. Johnson\thanks{Supported in part by NSF DMS-1301604  and the U.S.-Israel Binational Science
Foundation
 }, and   Gideon Schechtman\thanks{Supported in part by the U.S.-Israel Binational Science Foundation. Participant, 
 NSF Workshop in Analysis and Probability, Texas A\&M University
 }} \maketitle

\begin{abstract}
We construct a 
Schauder basis for $L_1$ consisting of non-negative functions and investigate unconditionally basic and quasibasic sequences of non-negative functions in $L_p$, $1\le p < \infty$.
\end{abstract}

\section{Introduction}
In \cite{ps}, Powell and Spaeth investigate non-negative sequences of functions in $L_p$, $1\le p < \infty$, that satisfy some kind of basis condition, with a view to determining whether such a sequence can span all of $L_p$.  They prove, for example, that there is no unconditional basis or even unconditional quasibasis (frame) for $L_p$ consisting of non-negative functions. On the other hand, they prove that there are non-negative  quasibases and non-negative $M$-bases  for $L_p$.  The most important question left open by their investigation is whether there is a (Schauder) basis for $L_p$ consisting of non-negative functions. In section \ref{section2}
 we show that there is basis for $L_1$ consisting of non-negative functions.
 
 In section \ref{section3} we discuss the structure of unconditionally basic non-negative normalized sequences in $L_p$, $1\le p < \infty$.  The main result is that such a sequence is equivalent to the unit vector basis of $\ell_p$. We also prove that the closed span  in $L_p$ of any unconditional quasibasic sequence embeds isomorphically into $\ell_p$.
 
 We use standard Banach space theory, as can be found in \cite{lt} or \cite{ak}.  Let us just mention that $L_p$ is $L_p(0,\infty)$, but inasmuch as this space is isometrically isomorphic under an order preserving  operator to $L_p(\mu)$ for any separable purely non-atomic measure $\mu$, our choice of $L(0,\infty)$ rather than e.g. $L_p(0,1)$ is a matter of convenience. Again as a matter of convenience, in the last part of Section \ref{section3} we revert to using $L_p(0,1)$ as a model for $L_p$. 
 
\section{A Schauder basis for $L_1(0,\infty)$ consisting of non-negative functions}\label{section2}

For $j=1,2,\dots$ let $\{h_{n,i}^j\}_{n=0, i=1}^{\infty\ \ \ 2^n}$ be the mean zero $L_1$ normalized Haar functions on the interval $(j-1,j)$.  That is, for $n=0,1,\dots$, $i=1,2,\dots,2^n$,
\[
h_{n,i}^j(t)=\left\{\begin{array}{cc}
              2^n & j-1+\frac{2i-2}{2^{n+1}}<t<j-1+\frac{2i-1}{2^{n+1}} \\
              -2^n & j-1+\frac{2i-1}{2^{n+1}} <t<j-1+\frac{2i}{2^{n+1}}
              \\
              0 & otherwise
            \end{array}\right.
 \]
 The system  $\{h_{n,i}^j\}_{n=0, i=1, j=1}^{\infty\ \ \ 2^n\ \ \infty}$, in any order which preserves the lexicographic order of $\{h_{n,i}^j\}_{n=0, i=1}^{\infty\ \ \ 2^n}$ for each $j$,  constitutes a basis for the subspace
 of $L_1(0,\infty)$ consisting of all functions whose restriction to each interval $(j-1,j)$ have mean zero.
 To simplify notation, for each $j$ we shall denote  by $\{h^j_i\}_{i=1}^\infty$ the system $\{h_{n,i}^j\}_{n=0, i=1}^{\infty\ \ \ 2^n}$ in its lexicographic order. We shall also denote by $\{h_i\}_{i=1}^\infty$ the union of the systems $\{h^j_i\}_{i=1}^\infty$, $j=1,2,\dots$, in any order that respects the individual orders of each of the $\{h^j_i\}_{i=1}^\infty$.

 Let $\pi$ be any permutation of the natural numbers and for each $i\in \N$ let $F_i$ be the two dimensional space spanned by $2{\bf 1}_{(\pi(i)-1,\pi(i))}+|h_i|$ and $h_i$.

 \begin{pr}\label{pr1}
 $\sum_{i=1}^\infty F_i$ is an FDD of $\overline{\rm span}^{L_1}\{F_i\}_{i=1}^\infty$.
 \end{pr}
 \pf
 The assertion will follow from the following inequality, which holds for all scalars $\{a_i\}_{i=1}^\infty$ and $\{b_i\}_{i=1}^\infty$,
 \begin{equation}\label{eq:fdd}
  \begin{array}{rl}
 \frac12\sum_{i=1}^\infty|a_i|+\frac18\|\sum_{i=1}^\infty b_ih_i\|&\le
 \|\sum_{i=1}^\infty a_i(2{\bf 1}_{(\pi(i)-1,\pi(i))}+|h_i|)+\sum_{i=1}^\infty b_ih_i\|\cr
 &\le
 3\sum_{i=1}^\infty|a_i|+\|\sum_{i=1}^\infty b_ih_i\|.
 \end{array}
 \end{equation}
 The right   inequality in (\ref{eq:fdd}) follows easily from the triangle inequality. As for the left   inequality, notice that the conditional expectation projection onto the closed span of $\{{\bf 1}_{(i-1,i)}\}_{i=1}^\infty$ is of norm one and the complementary projection, onto the closed span of $\{h_i\}_{i=1}^\infty$, is of norm 2. It follows that
 \[
 \|\sum_{i=1}^\infty a_i(2{\bf 1}_{(\pi(i)-1,\pi(i))})+\sum_{i=1}^\infty b_ih_i\|\ge \max\{2\sum_{i=1}^\infty|a_i|,\frac12 \|\sum_{i=1}^\infty b_ih_i\|\}.
 \]
 Since $\|\sum_{i=1}^\infty a_i|h_i|\|\le \sum_{i=1}^\infty |a_i|$, we get
 \[
 \|\sum_{i=1}^\infty a_i(2{\bf 1}_{(\pi(i)-1,\pi(i))}+|h_i|)+\sum_{i=1}^\infty b_ih_i\|\ge
 \max\{\sum_{i=1}^\infty|a_i|,\frac14 \|\sum_{i=1}^\infty b_ih_i\|\}
 \]
 from which the left hand side inequality in (\ref{eq:fdd}) follows easily.
 \endpf

  \begin{pr}\label{pr2}
 Let $\pi$ be any permutation of the natural numbers and for each $i\in \N$ let $F_i$ be the two dimensional space spanned by $2{\bf 1}_{(\pi(i)-1,\pi(i))}+|h_i|$ and $h_i$. Then $\overline{\rm span}^{L_1}\{F_i\}_{i=1}^\infty$ admits a basis consisting of non-negative functions.
 \end{pr}
 \pf
 In view of Proposition \ref{pr1} it is enough to show that each $F_i$ has a two term basis consisting of non-negative functions and with uniform basis constant.
 Put $x_i=2{\bf 1}_{(\pi(i)-1,\pi(i))}+|h_i| + h_i$ and $y_i=2{\bf 1}_{(\pi(i)-1,\pi(i))}+|h_i| - h_i$. Then clearly $x_i,y_i\ge 0$ everywhere and $\|x_i\|=\|y_i\|=3$. We now distinguish two cases: If ${\bf 1}_{(\pi(i)-1,\pi(i))}$ is disjoint from  the support of $h_i$ then, for all scalars $a,b$,
 \[
 \|ax_i+by_i\|\ge\|a(|h_i|+h_i)+b(|h_i|-h_i)\|=|a|+|b|.
 \]
 If the support of $h_i$ is included in $(\pi(i)-1,\pi(i))$, Let $2^{-s}$ be the size of that support, $s\ge 0$. Then  for all scalars $a,b$,
 \begin{equation*}
 \begin{array}{rl}
 \|ax_i+by_i\|&\ge\|a(|h_i|+h_i)+b(|h_i|-h_i)+2(a+b){\bf 1}_{{\rm supp}(h_i)}\|\\
 &=2^{-s-1}(|(2^{s+1}+2)a+2b|+|(2^{s+1}+2)b+2a|\ge \max\{|a|,|b|\}.
 \end{array}
 \end{equation*}
 \endpf

 \begin{thm}
 $L_1(0,\infty)$, and consequently any separable $L_1$ space,  admits a Schauder basis consisting of non-negative functions.
 \end{thm}

 \pf
 When choosing the order on $\{h_i\}$ we can and shall assume that $h_1=h_{0,1}^1$;  i.e., the first mean zero Haar function on the interval $(0,1)$. Let $\pi$ be any permutation of $\N$ such that $\pi(1)=1$ and for $i>1$, if $h_i=h_{n,k}^j$ for some $n,k,$ and $j$ then $\pi(i)>j$.
 It follows that except for $i=1$ the support of $h_i$ is disjoint from the interval $(\pi(i)-1,\pi(i))$. It is easy to see that such a permutation exists.
 We shall show that under these assumptions $\sum_{i=1}^\infty F_i$ spans $L_1(0,\infty)$ and, in view of Proposition \ref{pr2},  this will prove the theorem for $L_1(0,\infty)$.
 First, since $\pi(1)=1$ we get that $3{\bf 1}_{(0,1)}=2{\bf 1}_{(\pi(1)-1,\pi(1))}+|h_1|\in F_1$, and  since all the mean zero Haar functions on $(0,1)$ are clearly in $\sum_{i=1}^\infty F_i$,  we get that $L_1(0,1)\subset \sum_{i=1}^\infty F_i$.

 Assume by induction that $L_1(0,j)\subset \sum_{i=1}^\infty F_i$. Let $l$ be such that $\pi(l)=j+1$. By our assumption on $\pi$, the support of $h_l$ is included in $(0,j)$, and so by the induction hypothesis,  $|h_l|\in \sum_{i=1}^\infty F_i$. Since also $2{\bf 1}_{(j,j+1)}+|h_l|\in \sum_{i=1}^\infty F_i$ we get that ${\bf 1}_{(j,j+1)}\in \sum_{i=1}^\infty F_i$. Since the mean zero Haar functions on $(j,j+1)$ are also in $\sum_{i=1}^\infty F_i$ we conclude  that $L_1(0,j+1)\subset \sum_{i=1}^\infty F_i$.

 This finishes the proof for $L_1(0,\infty)$. Since every separable $L_1$ space is order isometric to one of the spaces $\ell_1^k$, $k=1,2,\dots$, $\ell_1, L_1(0,\infty)$, $L_1(0,\infty)\bigoplus_1\ell_1^k$, $k=1,2,\dots,$
 or $L_1(0,\infty)\bigoplus_1\ell_1$,  and since the discrete $L_1$ spaces $\ell_1^k$, $k=1,2,\dots,$ and $\ell_1$ clearly have non-negative bases, we get the conclusion for any separable $L_1$ space.
 \endpf
 
 \section{Unconditional non-negative sequences in $L_p$}\label{section3}
 Here we prove
 
  \begin{thm}\label{unconditional}
  Suppose that $\{x_n\}_{n=1}^\infty$ is a normalized unconditionally basic sequence of non-negative functions in $L_p$, $1\le p < \infty$.  Then $\{x_n\}_{n=1}^\infty$ is equivalent to the unit vector basis of $\ell_p$.
  \end{thm}
 \pf First we give a sketch of the proof, which should be enough for experts in Banach space theory. 
 By unconditionality, we have for all coefficients ${a_n}$ that $\| \sum_n a_n x_n \|_p$ is equivalent to the square function $\|(\sum_n |a_n|^2x_n^2)^{1/2}\|_p$, and, by non-negativity of $x_n$, is also equivalent to $\|\sum_n |a_n| x_n\|_p$. Thus by trivial interpolation when $1 \le p \le 2$, and by extrapolation when $2< p<\infty$, we see that 
 $\|\ \sum_n a_n x_n \|_p$ is equivalent to $\|(\sum_n |a_n|^p x_n^p)^{1/p}\|_p = (\sum_n |a_n|^p)^{1/p}$.
 
 We now give a formal argument for the benefit of readers who are not familiar with the background we assumed when giving the sketch.
 Let $K$ be the unconditional constant of $\{x_n\}_{n=1}^\infty$. Then
 \begin{equation}\label{eq3.1}
 \begin{split}
 K^{-1}\| \sum_{n=1}^N a_n x_n \|_p  \le B_p \| (\sum_{n=1}^N |a_n|^2x_n^2)^{1/2} \|_p \\  \le  B_p  \| \sum_{n=1}^N |a_n| x_n \|_p 
  \le B_p K \| \sum_{n=1}^N a_n x_n \|_p,
  \end{split}
 \end{equation}
 where the first inequality is obtained by integrating against the Rademacher functions (see, e.g., \cite[Theorem 2.b.3]{lt}).  The constant $B_p$ is Khintchine's constant, so $B_p=1$ for $p\le 2$ and $B_p$ is of order $\sqrt{p}$ for $p>2$. 
 If $1 \le p \le 2$  we  get from (\ref{eq3.1})   
  \begin{equation}\label{eq3.2}
K^{-1}\| \sum_{n=1}^N a_n x_n \|_p  \le \| (\sum_{n=1}^N |a_n|^p x_n^p)^{1/p} \|_p \le K \| \sum_{n=1}^N a_n x_n \|_p. 
 \end{equation}
Since $ \| (\sum_{n=1}^N |a_n|^p x_n^p)^{1/p} \|_p = (\sum_{n=1}^N |a_n|^p)^{1/p}$, this completes the proof when $1\le p \le 2$.  
 When $2<p<\infty$, we need to extrapolate rather than do (trivial) interpolation. Write $1/2 = \theta/1 + (1-\theta)/p$. Then  
  \begin{equation}\label{eq3.3}
  \begin{split}
(KB_p)^{-1 }  & \|    \sum_{n=1}^N a_n x_n \|_p   \le   \| (\sum_{n=1}^N |a_n|^2x_n^2)^{1/2} \|_p \\ & \le  \| \sum_{n=1}^N |a_n| x_n \|_p^\theta  \| (\sum_{n=1}^N |a_n|^p x_n^p)^{1/p} \|_p^{1-\theta} \\
& \le K \|    \sum_{n=1}^N a_n x_n \|_p^\theta     (\sum_{n=1}^N |a_n|^p)^{(1-\theta)/p}        , \quad  \text{so that} \\
(K^2 B_p)^{(-1)/(1-\theta)}  & \|    \sum_{n=1}^N a_n x_n \|_p    \le  (\sum_{n=1}^N |a_n|^p)^{1/p} \le K  \|    \sum_{n=1}^N a_n x_n \|_p. \quad \endpf
\end{split} 
 \end{equation}

 \bigskip
 
  As stated, Theorem \ref{unconditional} gives no information when $p=2$ because every normalized unconditionally basic sequence in a Hilbert space is equivalent to the unit vector basis of $\ell_2$. However, if we extrapolate slightly differently in the above argument (writing $1/2 = \theta/1 + (1-\theta)/\infty$) we see that, no matter what $p$ is, $\|    \sum_{n=1}^N a_n x_n \|_p$ is also equivalent to $\|\max_n |a_n| x_n \|_p$. From this one can deduce e.g. that only finitely many Rademachers can be in the closed span of $\{x_n\}_{n=1}^\infty$; in particular, $\{x_n\}_{n=1}^\infty$ cannot be a basis for $L_p$ even when $p=2$.  However, the proof given in    \cite{ps} that a normalized unconditionally basic sequence of non-negative functions $\{x_n\}_{n=1}^\infty$ in $L_p$ cannot span $L_p$ actually shows that   only finitely many Rademachers can be in the closed span of $\{x_n\}_{n=1}^\infty$.  This is improved in our last result, which shows that the closed span of an unconditionally  non-negative quasibasic sequence in $L_p(0,1)$ cannot contain any strongly embedded infinite dimensional subspace (a subspace $X$ of $L_p(0,1)$ is said to be strongly embedded if the $L_p(0,1)$ norm is equivalent to the $L_r(0,1)$ norm on $X$ for some -- or, equivalently, for all -- $r<p$; see e.g. \cite[p. 151]{ak}). The main work for proving this is   contained in  Lemma \ref{lemma3.1}.  
  
  Before stating Lemma \ref{lemma3.1}, we recall that a quasibasis for a Banach space $X$ is a sequence 
  $\{f_n,g_n\}_{n=1}^\infty$ in $X \times X^*$ such that for each $x$ in $X$ the series $\sum_n \langle g_n, x \rangle f_n$ converges to $x$. (In \cite{ps}  a sequence $\{f_n\}_{n=1}^\infty$ in $X$ is a called a quasibasis for $X$ provided there exists such a sequence $\{g_n\}_{n=1}^\infty$.  Since the sequence $\{g_n\}_{n=1}^\infty$ is typically not unique, we prefer to specify it up front.)  The quasibasis  $\{f_n,g_n\}_{n=1}^\infty$ is said to be unconditional provided that for each $x$ in $X$  the series $\sum_n \langle g_n, x \rangle f_n$ converges unconditionally to $x$. One then gets from the uniform boundedness principle (see, e.g., \cite[Lemma 3.2]{ps}) that there is a constant $K$ so that for all $x$ and all scalars $a_n$ with $|a_n|\le 1$, we have $\| \sum_n a_n  \langle g_n, x \rangle f_n \| \le K \|x\|$.   A sequence $\{f_n,g_n\}_{n=1}^\infty$ in $X \times X^*$  is said to be  [unconditionally] quasibasic provided $\{f_n,h_n\}_{n=1}^\infty$ is an [unconditional] quasibasis for the closed span $[f_n]$ of $\{f_n\}_{n=1}^\infty$, where $h_n$ is the restriction of $g_n$ to $[f_n]$.
  \begin{lm}\label{lemma3.1}
 Suppose that $\{f_n,g_n\}_{n=1}^\infty$ is an uncondtionally quasibasic sequence in $L_p(0,1)$, $1<p<\infty$ with each $f_n$ non-negative.  If $\{y_n\}_{n=1}^\infty$ is a normalized weakly null sequence in $[f_n]$, then $\|y_n\|_1 \to 0$ as $n\to \infty$.
  \end{lm}
  \pf
  If the conclusion is false, we get a normalized  weakly null sequence $\{y_n\}_{n=1}^\infty$ in $[f_n]$ and a $c>0$ so that for all $n$ we have $\|y_n\|_1 > c$.

  By passing to a subsequence of $\{y_n\}_{n=1}^\infty$, we can assume that there are integers $0=m_1 < m_2 < \dots$ so that for each $n$,
  \begin{equation}\label{yn}
    \sum_{k=1}^{m_n} |\langle g_k, y_n \rangle| \|f_k\|_p <  2^{-n-3}c \quad \text{and} \quad 
  \|\sum_{k={m_{n+1}+1}}^\infty|\langle g_k, y_n \rangle| f_k \|_p <   2^{-n-3}c.
  \end{equation}
  Effecting the first inequality in (\ref{yn}) is no problem because $y_n \to 0$ weakly, but the second inequality perhaps requires a comment.  Once we have a  $y_n$ that satisfies the   first inquality in (\ref{yn}), from the unconditional convergence of the expansion of $y_n$ and the non-negativity of all $f_k$ we get that 
  $\|\sum_{k=N}^\infty |\langle g_k, y_n \rangle| f_k \|_p \to 0$ as $n\to \infty$, which allows us to select $m_{n+1}$ to satisfy the second inequality in (\ref{yn}).

  Since $\|y_n\|_1 > c$, from (\ref{yn}) we also have for every $n$ that
  \begin{equation}\label{yn2}
   \| \sum_{k=m_n+1}^{m_{n+1}} |\langle g_k, y_n \rangle| f_k \|_1 \ge 
 \| \sum_{k=m_n+1}^{m_{n+1}} \langle g_k, y_n \rangle f_k \|_1 \ge c/2.
  \end{equation}

  Since $L_p$ has an unconditional basis,  by passing to a further subsequence we can assume that $\{y_n\}_{n=1}^\infty$  is unconditionally basic with constant $K_p$. Also, $L_p$ has type $s$, where $s = p \wedge 2$ (see \cite[Theorem 6.2.14]{ak}),  so 
  for some constant $K'_p$ we have for every $N$ the inequality 
  \begin{equation}\label{type}
  \| \sum_{n=1}^N y_n \|_p \le K'_p N^{1/s}. 
  \end{equation}
  On the other hand, letting $\delta_k = \text{sign}\,  \langle g_k, y_n \rangle $ when $m_n+1 \le k \le m_{n+1}$, $n =1,2,3,\dots$, we have
 \begin{equation}\label{lowerestimate}
   \begin{split}
K K_p \| \sum_{n=1}^N y_n \|_p   & \ge    K_p \| \sum_{n=1}^N   \sum_{k=1}^\infty \delta_k  \langle g_k, y_n \rangle f_k \|_p 
\\
& \ge \| \sum_{n=1}^N \sum_{k=m_n+1}^{m_{n+1}} |\langle g_k, y_n \rangle | f_k \|_p
- \| \sum_{n=1}^N \sum_{k\not\in [m_n+1, m_{n+1}]} \delta_k  \langle g_k, y_n \rangle  f_k \|_p
\\
& \ge 
 \| \sum_{n=1}^N \sum_{k=m_n+1}^{m_{n+1}} |\langle g_k, y_n \rangle | f_k \|_1 
 -  \| \sum_{n=1}^N \sum_{k\not\in [m_n+1, m_{n+1}]}  |  \langle g_k, y_n \rangle  | f_k \|_p
\\
& \ge 
  \sum_{n=1}^N \| \sum_{k=m_n+1}^{m_{n+1}} |\langle g_k, y_n \rangle | f_k \|_1 
  \\
 & - \sum_{n=1}^N \left( \sum_{k=1}^{m_n} | \langle g_k, y_n \rangle | \|f_k\|_p +   \|\sum_{k=m_{n+1}+1}^\infty | \langle g_k, y_n \rangle | f_n \|_p \right)
\\
&\ge 
Nc/2 - c/4 \quad\quad \text{by (\ref{yn2}) and (\ref{yn})}
\end{split}
  \end{equation}
This contradicts (\ref{type}).  \endpf

\begin{thm}\label{quasibasic}
Suppose that $\{f_n,g_n\}_{n=1}^\infty$ is an unconditional quasibasic sequence in $L_p(0,1)$, $1\le p < \infty$,  and each $f_n$ is non-negative.  Then the closed span $[f_n]$ of  $\{f_n\}_{n=1}^\infty$ embeds isomorphically into $\ell_p$.
\end{thm}
\pf
The case $p=1$ is especially easy:  Assume, as we may, that  $\|f_n\|_1=1$. There is a constant $K$ so that for each $y$ in $[f_n]$
\begin{equation}
\|y\|_1 \le \| \sum_{n=1}^\infty | \langle g_n, y \rangle | f_n \|_1 \le K \|y\|_1,
\end{equation}
hence the mapping $y \mapsto \{ \langle g_k , y \rangle \}_{k=1}^\infty$ is an isomorphism from $[f_n]$ into $\ell_1$.

So in the sequel assume that $p>1$.  From Lemma \ref{lemma3.1} and standard arguments (see, e.g., \cite[Theorem 6.4.7]{ak}) we have that every normalized weakly null sequence in $[f_n]$ has  a subsequence 
that is an arbitrarliy small perturbation of a disjoint sequence and hence the subsequence is $1+\epsilon$-equivalent to the unit vector basis for $\ell_p$. This implies that $[f_n]$ embeds isomorphically into $\ell_p$ (see \cite{jo} for the case $p>2$ and \cite[Theorems III.9, III.1, and III.2]{j} for the case $p<2$). \endpf


 \begin{tabular}{l}
W.~B.~Johnson\\
Department of Mathematics\\
Texas A\&M University\\
College Station, TX  77843 U.S.A.\\
{\tt johnson@math.tamu.edu}
\\
\end{tabular}

\bigskip

\begin{tabular}{l}
G.~Schechtman\\
Department of Mathematics\\
Weizmann Institute of Science\\
Rehovot, Israel\\
{\tt gideon@weizmann.ac.il}\\
\end{tabular}


\begin{thebibliography}{999}
  
  \bibitem{ak}
  Albiac, Fernando; Kalton, Nigel J. 
  Topics in Banach space theory. 
  Graduate Texts in Mathematics, 233. 
  Springer, New York, 2006.
  
   \bibitem{j}
  Johnson, W. B. 
  On quotients of $L_p$  which are quotients of $\ell_p$. 
  Compositio Math. 34 (1977), no. 1, 69Ð89.
  
  
   \bibitem{jo}
   Johnson, W. B.; Odell, E. 
   Subspaces of $L_p$ which embed into $\ell_p$. 
   Compositio Math. 28 (1974), 37Ð49.

\bibitem{lt}
Lindenstrauss, Joram; Tzafriri, Lior.
 Classical Banach spaces. I. Sequence spaces. 
 Ergebnisse der Mathematik 
  und ihrer Grenzgebiete, Vol. 92. 
  Springer-Verlag, Berlin-New York, 1977.



\bibitem{ps}
Powell, Alexander M.:    Spaeth, Anneliese H. 
Nonnegative constraints for spanning systems. Trans. AMS (to appear).



 \end{thebibliography}
\end{document}